\documentclass[12pt]{amsart}
\usepackage{amssymb,amsthm,amsmath,epsfig,latexsym}
\usepackage{calc}

\usepackage{ifthen}
\newcommand{\showcomments}{yes}
\newsavebox{\commentbox}
\newenvironment{comment}%
{\ifthenelse{\equal{\showcomments}{yes}}%
{\footnotemark
    \begin{lrbox}{\commentbox}
    \begin{minipage}[t]{1.25in}\raggedright\sffamily\upshape\tiny
    \footnotemark[\arabic{footnote}]}
{\begin{lrbox}{\commentbox}}}%
{\ifthenelse{\equal{\showcomments}{yes}}%
{\end{minipage}\end{lrbox}\marginpar{\usebox{\commentbox}}}
{\end{lrbox}}}

\begin{document}

\newcommand{\mmbox}[1]{\mbox{${#1}$}}
\newcommand{\proj}[1]{\mmbox{{\mathbb P}^{#1}}}
\newcommand{\Cr}{C^r(\Delta)}
\newcommand{\CR}{C^r(\hat\Delta)}
\newcommand{\affine}[1]{\mmbox{{\mathbb A}^{#1}}}
\newcommand{\Ann}[1]{\mmbox{{\rm Ann}({#1})}}
\newcommand{\caps}[3]{\mmbox{{#1}_{#2} \cap \ldots \cap {#1}_{#3}}}
\newcommand{\N}{{\mathbb N}}
\newcommand{\Z}{{\mathbb Z}}
\newcommand{\R}{{\mathbb R}}
\newcommand{\Tor}{\mathop{\rm Tor}\nolimits}
\newcommand{\Ext}{\mathop{\rm Ext}\nolimits}
\newcommand{\Hom}{\mathop{\rm Hom}\nolimits}
\newcommand{\im}{\mathop{\rm Im}\nolimits}
\newcommand{\rank}{\mathop{\rm rank}\nolimits}
\newcommand{\supp}{\mathop{\rm supp}\nolimits}
\newcommand{\arrow}[1]{\stackrel{#1}{\longrightarrow}}
\newcommand{\CB}{Cayley-Bacharach}
\newcommand{\coker}{\mathop{\rm coker}\nolimits}
\sloppy
\newtheorem{defn0}{Definition}[section]
\newtheorem{prop0}[defn0]{Proposition}
\newtheorem{conj0}[defn0]{Conjecture}
\newtheorem{thm0}[defn0]{Theorem}
\newtheorem{lem0}[defn0]{Lemma}
\newtheorem{corollary0}[defn0]{Corollary}
\newtheorem{example0}[defn0]{Example}

\newenvironment{defn}{\begin{defn0}}{\end{defn0}}
\newenvironment{prop}{\begin{prop0}}{\end{prop0}}
\newenvironment{conj}{\begin{conj0}}{\end{conj0}}
\newenvironment{thm}{\begin{thm0}}{\end{thm0}}
\newenvironment{lem}{\begin{lem0}}{\end{lem0}}
\newenvironment{cor}{\begin{corollary0}}{\end{corollary0}}
\newenvironment{exm}{\begin{example0}\rm}{\end{example0}}

\newcommand{\defref}[1]{Definition~\ref{#1}}
\newcommand{\propref}[1]{Proposition~\ref{#1}}
\newcommand{\thmref}[1]{Theorem~\ref{#1}}
\newcommand{\lemref}[1]{Lemma~\ref{#1}}
\newcommand{\corref}[1]{Corollary~\ref{#1}}
\newcommand{\exref}[1]{Example~\ref{#1}}
\newcommand{\secref}[1]{Section~\ref{#1}}

\newcommand{\std}{Gr\"{o}bner}
\newcommand{\jq}{J_{Q}}



\title {Smooth planar $r$-splines of degree $2r$}

\author{\c Stefan O. Toh\v aneanu}
\address{Mathematics Department \\ Texas A\&M University \\
  College Station \\ TX 77843-3368 \\ USA}
\email{tohanean@math.tamu.edu}


\begin{abstract}
\noindent In \cite{as}, Alfeld and Schumaker give a formula for
the dimension of the space of piecewise polynomial functions
(splines) of degree $d$ and smoothness $r$ on a generic
triangulation of a planar simplicial complex $\Delta$ (for $d \ge
3r+1$) and any triangulation (for $d\geq 3r+2$). In \cite{ss}, it
was conjectured that the Alfeld-Schumaker formula actually holds
for all $d \ge 2r+1$. In this note, we show that this is the best
result possible; in particular, there exists a simplicial complex
$\Delta$ such that for any $r$, the dimension of the spline space
in degree $d=2r$ is not given by the formula of \cite{as}. The
proof relies on the explicit computation of the nonvanishing of
the first local cohomology module described in \cite{ss2}.
\end{abstract}
\maketitle


\section{Introduction and preliminaries}
Let $\Delta$ be a planar, strongly-connected finite simplicial
 complex. The set of piecewise polynomial functions on $\Delta$
of smoothness $r$ has the structure of a module $\Cr$ over the
polynomial ring $\mathbb{R}[x,y]$; the subset of $f \in \Cr$ such
that $f|_{\sigma}$ is of degree at most $d$ for each two-simplex
$\sigma \in \Delta$ is a finite dimensional real vector-space,
denoted $\Cr_d$. In \cite{as}, for almost all triangulations,
Alfeld and Schumaker give a formula for the dimension of the
$\Cr_d$ in terms of combinatorial and local geometric data (data
depending only on local geometry at the interior vertices of
$\Delta$), as long as $d \ge 3r+1$. In \cite{ss}, it was
conjectured that their formula actually
 holds as long as $d \ge 2r+1$.The purpose of this brief note
is to show that this conjecture is optimal; in particular we
exhibit a planar simplicial complex $\Delta$ such that for
{\it any} $r$ there exist ``special'' splines
in degree $d=2r$. In other words, the conjecture of \cite{ss} is
tight.

The methods we use depend on the homological approach developed by
Billera in \cite{b} to answer a conjecture of Strang. This
approach was further developed by Schenck and Stillman in the
papers \cite{ss1}, \cite{ss2}, using a chain complex different
from Billera's and some additonal technical tools (local
cohomology and duality). When $\Delta$ is a planar simplicial
complex, it turns out that the delicate geometry of the problem is
captured by a certain local cohomology module, which as shown in
\cite{ss2} has a simple description (see below). For example, the
geometry of the famous ``Morgan-Scott'' example, which shows that
even for $r=1$ the formula of \cite{as} does not apply for $d=2$,
is captured by this local cohomology. For more results on the
Morgan-Scott triangulation see \cite{d}, \cite{dfk}, and for
computations of the dimension for small triangulations see
\cite{a}. In the next section we quickly review the presentation
of this module; then we exhibit a specific $\Delta$ and prove that
for any $r$, the dimension of this module (which exactly captures
the discrepancy between Alfeld-Schumaker's formula and the actual
dimension) is nonzero in degree $d=2r$.

\section{Review of local cohomology}
By taking the cone $\hat \Delta$ over $\Delta$, we turn the problem
of computing $\dim \Cr_d$ into a problem in commutative
algebra - compute the Hilbert function $\CR_d$. As shown in
\cite{ss2}, $$\dim \Cr_d = \dim \CR_d=L(\Delta,r,d)+\dim N_d$$ where $L(\Delta,r,d)$
is the Alfeld-Schumaker formula and $N$ is a graded
$R=\mathbb{R}[x,y,z]$ module of finite length. Lemma 3.8 of \cite{ss2}
contains the following description: $N$ is the quotient of a free
module generated by the totally interior edges (those edges with no
vertex $ \subseteq \partial \Delta$), modulo the syzygies at each
interior vertex. The generators of $N$ are shifted so that they have
degree $r+1$. This description seems cumbersome, but as we'll see
in the example below, it is fairly easy to work with.

So in the terms above, the conjecture of \cite{ss} is that $N$
 vanishes in degree $2r+1$. Our goal is to show that this bound
is the best possible, so we want to find a configuration $\Delta$
such that for all $r$, $N_{2r} \ne 0$. Consider the
following simplicial complex (see also \cite{ss})
\begin{center}
\epsfig{figure=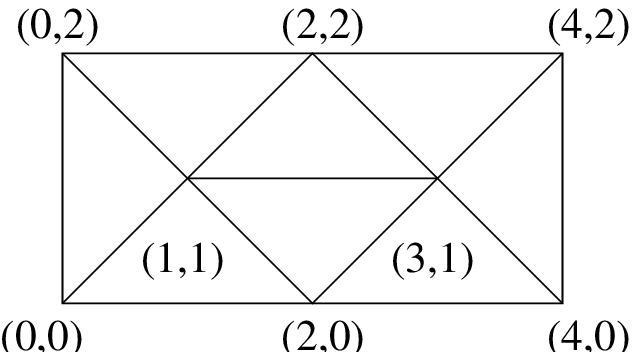,height=1in,width=1.7in}
\end{center}

To find $N$, we begin by determining the minimal free resolutions
for the ideals $I_i=\mathcal{J}(v_i)$ for $v_1$ and $v_2$ the
interior vertices. We have:
$$I_1=\langle (x+y-2z)^{r+1},(x-y)^{r+1},(y-z)^{r+1}\rangle$$
$$I_2=\langle (x+y-4z)^{r+1},(x-y-2z)^{r+1},(y-z)^{r+1}\rangle.$$
These ideals are in $R=\mathbb{R}[x,y,z]$. Notice that $y-z$ is the
linear form vanishing on the totally interior edge.
With the change of variables given by the matrix
$$
\left[
\begin{array}{rrr}
1& 1& -2 \\
0& -2& 2 \\
1& 1& -4
\end{array}
\right]
$$
we can suppose that
$$I_1=\langle x^{r+1},(x+y)^{r+1},y^{r+1}\rangle, \mbox{ and }
I_2=\langle z^{r+1},(z+y)^{r+1},y^{r+1}\rangle.$$

The minimal free resolutions for these ideals are:
$$0 \longrightarrow R^2 \stackrel{\tiny{
\left[
\begin{array}{c}
A_1\, D_1 \\
B_1\, E_1 \\
C_1\, F_1
\end{array}
\right]}
}
{\longrightarrow} R^3  \stackrel{\tiny{
\left[
\begin{array}{c}
x^{r+1}\, (x+y)^{r+1}\, y^{r+1}
\end{array}
\right]}
}
{\longrightarrow} R \longrightarrow R/I_1 \longrightarrow 0$$ and
$$0 \longrightarrow R^2 \stackrel{\tiny{
\left[
\begin{array}{c}
A_2\, D_2 \\
B_2\, E_2 \\
C_2\, F_2
\end{array}
\right]}
}
{\longrightarrow} R^3  \stackrel{\tiny{
\left[
\begin{array}{c}
z^{r+1}\, (z+y)^{r+1}\, y^{r+1}
\end{array}
\right]}
}
{\longrightarrow} R \longrightarrow R/I_2 \longrightarrow 0.$$

By \cite{ss2}, Lemma 3.8, $N\approx R(-r-1)/\langle
C_1,F_1,C_2,F_2\rangle$. In what follows we will prove that the
Hilbert function of $R/\langle C_1,F_1,C_2,F_2\rangle$ is nonzero
in degree $r-1$, for any positive integer $r$. In other words
\begin{eqnarray}
HF(N,2r)
&=& HF(R(-r-1)/\langle C_1,F_1,C_2,F_2\rangle,2r)\nonumber \\
&=& HF(R/\langle C_1,F_1,C_2,F_2\rangle,r-1)\nonumber \\
&\ne& 0.\nonumber
\end{eqnarray}

\section{The Hilbert function of $R/\langle C_1,F_1,C_2,F_2\rangle$ is
nonzero in degree $r-1$}

In the previous section we saw that the ideals $I_1$ and $I_2$
have a special form. First notice that they are symmetric (in
terms of the generators) in $x$ and $z$. So replacing $x$ by $z$
in the forms of $C_1$ and $F_1$ we obtain $C_2$ and $F_2$. Next
observe that we can look at the ideal $I_1$ as an ideal in
$A=\mathbb{R}[x,y]$. Similarly, $I_2$ is an ideal
$A'=\mathbb{R}[y,z]$. Hence $C_1,F_1 \in A$ and $C_2,F_2 \in A'$.

For $i=1,2$, the ideal $\langle C_i,F_i \rangle$ is a complete
intersection. For example, if $\langle C_1,F_1 \rangle$ is not a
complete intersection, since $C_1 \neq 0$ and $F_1 \neq 0$, then
$codim(\langle C_1,F_1 \rangle)=1$. Therefore, there is a nonunit
in $A$, say $d_1$, such that $d_1 | C_1$ and $d_1 | F_1$.
Therefore, $d_1 | x^{r+1}=B_1F_1-E_1C_1$ and $d_1 |
(x+y)^{r+1}=A_1F_1-D_1C_1$. Hence $codim(\langle
x^{r+1},(x+y)^{r+1}\rangle)=1$. But this contradicts the fact that
$codim(\langle x^{r+1},(x+y)^{r+1}\rangle)=2$, as
$\{x^{r+1},(x+y)^{r+1}\}$ is a regular $A$-sequence. So the ideal
$\langle C_1,F_1 \rangle$ is a complete intersection. The same
argument shows that $\langle C_2,F_2 \rangle$ is also a complete
intersection.

These observations will simplify our future computations. We'll
need to discuss two cases, depending on if $r$ is odd or even.

\subsection{$r+1=2n$}
Let $A=\mathbb{R}[x,y]$ and
$I_1=\langle x^{2n},(x+y)^{2n},y^{2n}\rangle$. By \cite{ss1}, Theorem 3.1, a free
resolution for $I_1$ is:
\vskip .2in
$ 0 \longrightarrow
A(-3n)^2 \stackrel{\tiny{
\left[
\begin{array}{c}
A_1\, D_1 \\
B_1\, E_1 \\
C_1\, F_1
\end{array}
\right]}
}
{\longrightarrow} A(-2n)^3  \stackrel{\tiny{
\left[
\begin{array}{c}
x^{2n}\, (x+y)^{2n}\, y^{2n}
\end{array}
\right]}
}
{\longrightarrow}  I_1
\longrightarrow 0$,
\vskip .2in
where
$\deg C_1=\deg F_1=3n-2n=n$. From the observations at the beginning of this
section we get the minimal free resolution for $A/\langle C_1,F_1\rangle$:
$$
0 \longrightarrow A(-2n) \longrightarrow A(-n)^2 \longrightarrow  A
\longrightarrow A/\langle C_1,F_1\rangle.$$
Therefore the Hilbert series is
$HS(A/\langle C_1,F_1\rangle,t)=\frac{1-2t^n+t^{2n}}{(1-t)^2}$. Hence there exists a
monomial $x^uy^v$ of degree $2n-2=r-1$ which is not in
$\langle C_1,F_1\rangle$. In fact, with an easy computation
we can see that this monomial is actually $x^{2n-2}$.

All we did here is in two variables $x$ and $y$. Let's go back to the ring
$R=\mathbb{R}[x,y,z]$ and suppose that the above monomial $x^uy^v$ is in
$\langle C_1,F_1\rangle+\langle C_2,F_2\rangle=\langle C_1,F_1,C_2,F_2\rangle$. Then there are
${\alpha}_1,{\beta}_1,{\alpha}_2,{\beta}_2 \in R$ such that:
$$x^uy^v={\alpha}_1C_1+{\beta}_1F_1+{\alpha}_2C_2+{\beta}_2F_2.$$
In this equation, since for $x=z$ we get $C_1=C_2$ and $F_1=F_2$
(see the remarks at the beginning), we obtain an equation in
$A=\mathbb{R}[x,y]$:
\begin{eqnarray}
x^uy^v &=&
{\alpha}_1(x,y,x)C_1+{\beta}_1(x,y,x)F_1+{\alpha}_2(x,y,x)C_1+{\beta}_2(x,y,x)F_1
\nonumber \\
&=& {\alpha}_1'C_1+{\beta}_1'F_1 \nonumber
\end{eqnarray}

So $x^uy^v \in \langle C_1,F_1\rangle$. This contradicts the way
we chose $x^uy^v$. Hence there is a monomial of degree $r-1$ which
is not in $\langle C_1,F_1,C_2,F_2\rangle$.

\subsection{$r+1=2n+1$}
For the odd case, the idea is almost identical.
Let $A=\mathbb{R}[x,y]$ and
$I_1=\langle x^{2n+1},(x+y)^{2n+1},y^{2n+1}\rangle$. Again, by \cite{ss1}, Theorem 3.1, a free
resolution for $A/I_1$ is:
 \vskip .2in
$ 0 \longrightarrow
\begin{array}{c}
A(-3n-1)\\
\oplus\\
A(-3n-2) \\
\end{array}
\stackrel{\tiny{
\left[
\begin{array}{c}
A_1\, D_1 \\
B_1\, E_1 \\
C_1\, F_1
\end{array}
\right]}
}
{\longrightarrow} A(-2n-1)^3  \stackrel{\tiny{
\left[I_1
\right]}
}
{\longrightarrow}  A \longrightarrow A/I_1
\longrightarrow 0$,
 \vskip .2in
 where $\deg C_1=3n+1-(2n+1)=n$ and $\deg F_1=3n+2-(2n+1)=n+1$. $\langle C_1,F_1\rangle$ is a complete intersection so the minimal free resolution for $A/\langle C_1,F_1\rangle$ is:
$$
0 \longrightarrow A(-2n-1) \longrightarrow A(-n)\oplus A(-n-1)
\longrightarrow  A \longrightarrow A/\langle C_1,F_1\rangle.$$
Therefore the Hilbert series is $HS(A/\langle
C_1,F_1\rangle,t)=\frac{1-t^n-t^{n+1}+t^{2n+1}}{(1-t)^2}$. Hence
there exists a monomial $x^uy^v$ of degree $2n-1=r-1$ which is not
in $\langle C_1,F_1\rangle$. As in 3.1., the same argument gives
us that in fact this monomial is not in $\langle
C_1,F_1,C_2,F_2\rangle$.
\\

In conclusion the Hilbert function of $R/\langle
C_1,F_1,C_2,F_2\rangle$ is nonzero in degree $r-1$. This is
exactly what we wanted to see.
\vskip .2in \noindent {\bf
Acknowledgments} I thank Hal Schenck for
 helpful comments and discussions and for introduction to
the subject. This research was partially supported by NSF grant
DMS-03-11142 and ATP grant 010366-0103. I thank two anonymous
referees for useful suggestions.

\renewcommand{\baselinestretch}{1.0}
\small\normalsize 

\bibliographystyle{amsalpha}

\end{document}